\documentclass[11pt]{amsart}
\allowdisplaybreaks[1]

\usepackage{amssymb}
\usepackage{amsmath} 
\usepackage{amscd} 
\usepackage{epstopdf}

\theoremstyle{plain}
\newtheorem{theorem}{Theorem}

\theoremstyle{definition}
\newtheorem{definition}{Definition}
\newtheorem{example}{Example}

\newcommand{\ts}{\hspace{0.5pt}}

\newcommand{\RR}{\mathbb{R}\ts}

\newcommand{\ZZ}{\mathbb{Z}}

  \newcommand{\vol}{\mbox {vol}}
 \newcommand{\CalA}{\mathcal{A}}

\newcommand{\card}{\mathrm{card}}
\newcommand{\dens}{\mathrm{dens}}

\begin{document}

\title[Deforming Meyer sets]
{Deforming Meyer sets}
\dedicatory{For our friend Ludwig Danzer, on the occasion of his 80th birthday.}

\author{Jeong-Yup Lee}
\address{KIAS, 207-43, Cheongnyangni 2-dong, 
\newline 
\hspace*{12pt}Dongdaemun-gu, Seoul 130-722, Korea}
\email{jylee@kias.re.kr}

\author{Robert V.\ Moody}
\address{Department of Mathematics and Statistics, 
University of Victoria, \newline
\hspace*{12pt}Victoria, BC, V8P 5C2, Canada}
\email{rmoody@uvic.ca}

\begin{abstract} A linear deformation of a Meyer set $M$ in $\RR^d$ is
the image of $M$ under
a group homomorphism of the group $[M]$ generated by $M$
into $\RR^d$. We provide a necessary and sufficient condition for
such a deformation to be a Meyer set.  In the case that the deformation 
is a Meyer set and the deformation is injective, the deformation is 
pure point diffractive if the orginal set $M$ is pure point diffractive.
 \end{abstract}

\maketitle

AMS Classification Codes 52C23 51P05

\section{Introduction}

Meyer sets play an important role in the study of long-range aperiodic order. They are very easy to define ($M \subset \RR^d$ is a Meyer set if
it is relatively dense and $M-M$ is uniformly discrete \cite{Lagarias,M}), are generally aperiodic, and have an amazing degree of internal long-range order: they always have a relatively dense set of Bragg peaks \cite{Nicu}. It is this last property that makes them so attractive in the study of quasicrystals, since the prominent existence of Bragg peaks is an essential feature of quasicrystals \footnote{For example, we read in 1991 Report of the Ad Interim Commission on Aperiodic Crystals, which is a Commission of the International
Union of Crystallographers, ``by `crystal' we mean any solid having an essentially discrete diffraction diagram, and by `aperiodic crystal' we mean any crystal in which three-dimensional lattice periodicity can be considered to be absent.'' See also
\cite{Shechtman}.}.  All model sets and all relatively dense subsets of model sets are Meyer sets. For this reason Meyer sets are ubiquitous in the theory of long-range aperiodic order.

Recently, attention has been turning to the deformations of tilings and point sets in order to see the effect of deformation on diffraction and various topological invariants, with the hope of better understanding the nature of long-range aperiodic order \cite{BL2, BD, Sadun}. In this short note we extend this by looking at the effect of  additive deformations on Meyer sets. 

Briefly, any Meyer set $M$ in $\RR^d$ generates a subgroup $[M]$ of finite rank
(see \cite{M}, Prop.~7.4), normally larger than $d$, and by an additive deformation of $M$ we mean
a $\ZZ$-homomorphism of $[M]$ into $\RR^d$. Such a homomorphism, even when it is injective, can map $M$ into a set which is in some bounded neighbourhood of a hyperplane in $\RR^d$ -- a set of the from $B + H$ where $B$ is a ball and $H$ a hyperplane, and the result is not relatively dense and hence not a Meyer set. However, when this does not happen, the resulting image is Meyer, whether or not the homomorphism is injective. Furthermore, assuming injectivity, if $M$ is pure point diffractive, then so is $f(M)$.

The key to the proof is a characterization of Meyer sets called the linear approximation property \cite{Meyer, M} (see below for the definition), which up to now has not been
used much in the study of Meyer sets.

\section{Linear deformation of Meyer sets}

 For a subset $M$ of $\RR^d$, we let $[M]$ be the subgroup of 
$\RR^d$ generated by $M$. If $M$ is Delone and has finite local complexity (FLC) (and, in particular, if $M$ is a Meyer set) then $[M] \simeq \ZZ^k$ for some $k\ge d$.

We say that a subset $S$ of $\RR^d$ is {\bf tied} to a hyperplane (i.e. to a $(d-1)$-dimensional subspace) $H$ of $\RR^d$ if
there is $r>0$ so that $S \subset H + B_r$. $S$ is {\bf untied} if not tied to any
hyperplane. 
If $T\subset V$ is a subset of some $\ZZ$-module $V$ in $\RR^d$ and $f:V \longrightarrow \RR^d$
is a $\ZZ$-homomorphism, then $f$ is called a {\bf tied map on $T$} if $f(T)$ is a tied subset of $\RR^d$.

\smallskip

\begin{example} \label{linear-deformation-example}
Consider the Meyer set \[M:= \{ x \in \ZZ[\tau] \,:\, x^\star \in [0,1] \} \, , \]
where $\tau := (1+\sqrt 5)/2$ and $(\phantom{n})^\star$ is the conjugation
map on $\ZZ[\tau]$ defined by $\tau \mapsto -1/\tau$. Here, since
$1,\tau \in M$, $[M] = \ZZ[\tau] \simeq \ZZ^2$
and the set of points
$\{(a,b) \in \ZZ^2 \,:\, a+ b\tau \in M \}$ is tied to the line $y= \tau x$.
The mapping
$(\phantom{n})^\star$ is an example of a tied map on $M$ since it sends $M$ into the closed ball of radius $1$ around the hyperplane $\{0\}$. In fact, an arbitrary  $\ZZ$-map
$f: \ZZ[\tau] \longrightarrow \RR$ is tied on $M$ if and only if $f(\tau)/f(1) = -1/\tau$. Thus $f$ is untied except when it is a multiple of the map $(\phantom{n})^\star$. Needless to say, the image of $M$ under a tied map is not a Meyer set. However, for any mapping of $[M]$ into $\RR$ which is not tied, e.g.,  $f: (a+ b\tau) \mapsto a \sqrt 2+ b \pi $, the image of $M$ is a Meyer set by Theorem \ref{fMMeyer} below.
\end{example}

\smallskip
\definition \cite{Meyer,M}  $M\subset \RR^d$ satisfies the {\em linear approximation property} 
if for each additive homomorphism into a linear space, $f: [M] \longrightarrow \RR^m$, there is a linear map 
\begin{eqnarray} \label{linear-map-associated}
F: \RR^d  \longrightarrow \RR^m
\end{eqnarray}
with the property that $|F(x) - f(x)|$ is uniformly bounded on $M$.  

\smallskip
It is clear that such a linear map $F$ must be unique. We call it the {\it linear map associated with} $f$.

\begin{theorem} \label{lap}\cite{M} 
A relatively dense subset $M \subset \RR^d$ is a Meyer set if and only if
$[M]$ is finitely generated and has the linear approximation property.
\qed
\end{theorem}

\begin{theorem} \label{fMMeyer}
Let $M \subset \RR^d$ be a Meyer set and let
\[ f:[M] \longrightarrow \RR^d \]
 be any $\ZZ$-homomorphism, i.e. a homomorphism of groups.
 Assume that $f(M)$ is untied. Then $f(M)$ is a Meyer set of $\RR^d$.
 \end{theorem}

\noindent
{\sc \bf Proof:} Using Thm.~1, let
 $F: \RR^d  \longrightarrow \RR^d$ be the linear map associated with $f$, and let $B\subset \RR^d$ be a closed ball centred on the origin so that for all $x \in M$, $ F(x) - f(x) \in B$. 
 
The mapping $F$ is onto, as otherwise $F(\RR^d) \subset H$ for some hyperplane $H$,
and then $f(M) \subset F(M) + B \subset H + B$ contradicts the assumption that $f$ is untied on $M$. Thus $F$ is an isomorphism of vector spaces and is a homeomorphism of topological spaces. 

Let $r>0$ be chosen so that the $r$-cube $C_r$ of side length $r$ centred at $0$ has
the property that for all $y \in \RR^d$, $(y+ C_r) \cap M \ne  \emptyset$. Let 
$y_i \in \RR^d$, $i=1,2, \dots$ be chosen so that 
\[ \bigcup_{i=1}^\infty (y_i + C_r) = \RR^d \]
be a tiling of $\RR^d$ by copies of $C_r$. Then 
\[ \bigcup_{i=1}^\infty  (F(y_i)  + F(C_r))  = \RR^d \] is a tiling of 
$\RR^d$ by translates of the parallelepiped $F(C_r)$ and
\[\bigcup_{i=1}^\infty  (F(y_i)  + F(C_r) + B) \]
 is a covering of $ \RR^d$ by translates of the compact set
$F(C_r) + B$. For each $i$ we know there is some $x \in (y_i + C_r) \cap M$,
and then $F(x) \in F(y_i) + F(C_r)$ and $f(x) \in (F(y_i) + F(C_r) + B)\cap f(M)$. This shows that $f(M)$ is relatively dense.

Next we show that $f(M)$ has finite local complexity, and as a consequence it also
is uniformly discrete. Let $K \subset \RR^d$  be compact. We wish to show that
$(f(M) - f(M)) \cap K$ is finite. All of its elements are  of the form $f(x) - f(y) \in K$. Then $F(x-y) = F(x) - F(y)  \in K + 2B$ so $x-y \in (M - M)  \cap F^{-1}(K + 2B)$, which is a
finite set. Thus $(f(M) - f(M)) \cap K \subset f((M - M)  \cap F^{-1}(K + 2B))$, which is a finite set, too. This proves finite local complexity and hence that $f(M)$ is a Delone set.

Since $M$ is Meyer, there is a finite set
$S$ so that $M-M \subset M + S$ \cite{Meyer, M}. We can clearly assume $S \subset [M]$. Then
$f(M) - f(M) \subset f(M) + f(S)$ shows that $f(M)$ is a Meyer set.
\qed

\smallskip
\noindent
{\sc \bf Remarks}
\begin{itemize}
\item[{1.}] In the last paragraph of the proof of Thm.~\ref{fMMeyer} we have used the 
well-known result of Lagarias that a Delone set $M \subset \RR^d$ is Meyer if and only if there is a finite set $S$ with $M-M \subset M+S$  \cite{Lagarias, M}.
\item[{2.}] Note that this result does not require $f$ to be $1$--$1$.
\item[{3.}] Under the hypotheses of Thm.~\ref{fMMeyer}, $M-M +M$ is also a Meyer set and it generates the same
subgroup of $\RR^d$ as $M$, namely $[M]$. Then $F$ is still the unique linear mapping associated with $f$. However, if $F(x) -f(x) \in B$ for all $x\in M$ then
$F(y) - f(y) \in 3B$ for all  $y \in M-M+M$. We will make use of the set $M- M + M$ in the next result.
\end{itemize}
 
 \smallskip
 
\begin{theorem} \label{pp}
Let $M \subset \RR^d$ be a Meyer set that is pure point diffractive 
with respect to some van Hove sequence 
 $\CalA = \{ A_m\}_{m=1}^\infty $. Let
\[ f:[M] \longrightarrow \RR^d \]
 be a $\ZZ$-homomorphism.
 Assume that $f$ is one-to-one on $M$ and that $f(M)$ is untied.  Let $F :[M] \longrightarrow \RR^d$ be the unique $\RR$-linear mapping associated with $f$ as in 
$\eqref{linear-map-associated}$.  
 Then $f(M)$ is pure point diffractive {\rm (}with respect to the averaging sequence $F\CalA = \{ F(A_m)\}_{m=1} ^\infty${\rm )}.
 \end{theorem}
 
\noindent
 {\sc \bf Proof:} We are assuming that $\CalA = \{ A_m\}_{m=1}^\infty$ is a van Hove sequence, that is, each $A_m$ is a measurable and pre-compact set and for all compact sets $K$,
 \[ \lim_{m\to\infty} \frac{\vol(\partial^KA_m)}{\vol(A_m)} = 0  \, ,\]
 where $\partial^KA_m$ is the $K$-boundary of $A_m$.\footnote{The $K$-boundary of a set $A$ is
$ (( K + A )\setminus A^\circ ) \cup (( - K +\overline{\RR^d\setminus A})\cap A )$.}
 Since $F$ is a homeomorphism, $F(\partial^K A_m) = \partial^{F(K)}F(A_m)$  and each $F(A_m)$ is pre-compact. Thus $F\CalA$ is also a van Hove sequence. 
 
 For $N \subset \RR^d$ we define 
 \[\dens_\CalA(N) = \lim_{n\to\infty} \frac{\card (N \cap A_m)}{\vol(A_m)} 
 \, ,\]
 if it exists.
  
 Now let $N \subset M - M + M$. Then
 \begin{alignat*}{2} \card (F(N) \cap F(A_m)) \le \card (f(N) &\cap &&(F(A_m) + 3B))\\
&\le &&\card (F(N) \cap (F(A_m) + 6B)) \, . 
 \end{alignat*}
 
 Dividing by $\vol \, F(A_m)$, which is $\vol(A_m)$ scaled by the Jacobian
 $|\det(F)|$ of $F$, and using the fact that $F$ and $f$
 are $1$--$1$ functions, we obtain
  \begin{eqnarray*}
\frac{ \card (N \cap A_m)}{\vol(A_m)} &=& 
|\det(F)|\left( \frac{\card (F(N) \cap F(A_m))}{\vol(F(A_m))}\right)\\  
&\le&  |\det(F)|\left( \frac{\card (f(N) \cap (F(A_m)+3B))}{\vol(F(A_m))}\right)\\
&\le& |\det(F)|\left( \frac{\card (F(N) \cap (F(A_m) + 6B))}{\vol(F(A_m))}\right) \, .
\end{eqnarray*}
 Using the van Hove property, the second and fourth terms are equal in the limit and
 we obtain
 \[ \dens_{F\CalA}( f(N)) = \frac{1}{|\det(F)|} \dens_{\CalA}(N)  \, .\]
 
Since $M$ is Meyer, to say that it is pure point diffractive (relative to $\CalA$) is equivalent to saying that
for all $\epsilon >0$
\[ P(\epsilon) := \{ t :  \dens_\CalA((t + M)\, \triangle\, M) < \epsilon\} \]
is relatively dense (including the statement that these densities exist) \cite{BM}.

Suppose that $\epsilon < 2\,\dens(M)$. Then the sets $P(\epsilon)$ are subsets of $M-M$.
Since $f$ is one-to-one, for $t \in P(\epsilon)$ we have 
\[x \in (t+M) \, \triangle \, M \Leftrightarrow f(x) \in f(t+M)\, \triangle \, f(M)\, .\] 
Thus
\[f((t+M)\, \triangle \, M) = (f(t) + f(M))\, \triangle\, f(M) \, .\]
Using the density result above, and writing $P_{f(M)}(\phantom{n})$ for the corresponding 
$P$-function on $f(M)$, we have
\[ P_{f(M)}(\epsilon/|\det F|) \supset f(P(\epsilon)) \,.\]
Since $f(P(\epsilon))+ 2B \supset F(P(\epsilon))$, which is relatively dense, so too
is $f(P(\epsilon))$ and then $P_{f(M)}(\epsilon/|\det (F)|)$ is relatively dense. 

This being true for all $\epsilon < 2\,\dens(M)$, it follows
from this that $f(M)$ is pure point diffractive relative to the new
averaging sequence \cite{BM}. \qed

\smallskip

In general, diffraction depends on the van Hove sequence over which the averaging is made.
Thus potentially  one needs to be aware of the two different averaging sequences that occur
in Theorem~\ref{pp}. In many situations it does not matter which averaging sequence is used. 
For example, a Meyer set $M$ has uniform patch frequencies if and only
if its dynamical hull is uniquely ergodic, a condition independent of averaging sequences,
see \cite{Schlottmann, LMS}. In this situation the autocorrelation and the diffraction are likewise 
independent of the van Hove averaging sequences.

Every Meyer set is a subset of a model set with
real internal space (see \cite{M}, Prop. 9.2) and hence has a  tied map.  Given the role of tied maps in the
deformation theory which we have just described, one might be led to suppose that the existence of
tied maps somehow implies the Meyer property. However, the following example shows that the existence of a tied map from a Delone point set, even one with FLC, does not imply the Meyer property of the set.
\begin{example} 
Consider the alphabetic substitution $a \rightarrow aba$ and $b \rightarrow aaaa$. Using relative
tile lengths given by the coefficients of the left Perron-Frobenius eigenvector of the corresponding substitution matrix, 
we can construct from it a substitution tiling $\mathcal{T}$ in $\RR$ with expansion map $\phi:\,\phi(x) = (1 + \sqrt{5})x$. Then $\mathcal{T}$ clearly has FLC, but it does not have the Meyer property, since $1+ \sqrt{5}$ is not a Pisot number (see \cite[Ch.2, Thm.~6]{Meyer}, \cite{Lagarias99,solomyak}).
We choose representative points for each of the two types of tiles in the tiling and derive from them
a non-Meyer set $M'$. 

Now we construct a point set $\widetilde{M} \subset \RR^2$: 
\[\widetilde{M} = \{(a, b) \in \RR^2 \,:\, a \in M', b \in M\} \]
where $M$ is the point set of Example \ref{linear-deformation-example} above. Then $\widetilde{M}$ still has FLC, but it does not have the Meyer property. Define the map 
$f : [\widetilde{M}] \longrightarrow \RR^2$ by $f(a, b) = (a, b^{\star})$, where $(\phantom{n})^{\star}$ is the conjugation map on $\ZZ[\tau]$ defined in 
Example \ref{linear-deformation-example}. Then $f$ is a tied map.
\end{example}

{\sc Acknowedgment} The authors would like to thank Michael Baake for his interest and valuable comments in the development of this paper.

\end{document}